\documentclass[12pt]{article}

\usepackage{amsmath,amssymb,amsfonts,setspace}
\usepackage{mathtools}

\usepackage{authblk}

\def\qed{\hfill\mbox{\rule{0.5em}{0.5em}}}

\newcommand{\eq}{\begin{equation}}
\newcommand{\en}{\end{equation}}

\newcommand{\prob}{{\mathbb P}}
\newcommand{\me}{{\mathbb E}}

\newtheorem{theorem}{Theorem}

\def \proof{\noindent{\it Proof.\ }}

\begin{document}

\author{Michael Farber, Alexander Gnedin and Wajid Mannan }
\affil{Queen Mary, University of London}
\title{A Random Graph Growth Model}
\date{}

\pagestyle{myheadings}
\maketitle

\begin{abstract}
\noindent
A  growing  random graph  is constructed by successively sampling without replacement an element from the pool of virtual vertices and edges.
At start of the process the pool contains $N$  virtual vertices and no edges. Each time a vertex is sampled and occupied, the edges linking the vertex to previously occupied
vertices are added to the pool of virtual elements. We focus on the edge-counting at times when the graph has $n\leq N$ occupied vertices.
Two different Poisson limits are identified for $n\asymp N^{1/3}$ and $N-n\asymp 1$. For the bulk of the process, when $n\asymp N$,
the scaled number of edges is shown to fluctuate about a deterministic curve, with fluctuations  being of the order of $N^{3/2}$ and  approximable by a Gaussian bridge.

\end{abstract}

\section{The graph model}

We consider a  random graph growing  according to simple rules.   The idea is to successively sample from a  variable   pool of  virtual vertices and edges.
The sampled elements, called in the sequel occupied, comprise a    graph.
Each time a new vertex is  sampled the pool is topped up with virtual   edges missing in the current graph. 
The method creates a bias in the probabilities of edges in a way  that the vertices  sampled earlier are more likely to get linked. Moreover, the probabilities of new    
edges depend on the content of   the pool.
These features distinguish our process from  numerous random graph models studied in the literature  \cite{RG, Durrett, Frieze, Janson, Remco,  Kolchin}.
Our interest to the new model is motivated by the analogous construction of random   simplicial complexes
 but we leave this extension for future work and in this note  only focus   on the  edge  counting.

In more detail, the process evolves as follows.
At start there is a  pool of $N$  virtual vertices and no edges.
At each subsequent step there is some number $n\leq N$ of occupied  vertices, which are linked by some  number  $k\leq {n\choose 2}$
of occupied  edges.      
The occupied vertices and edges comprise a   graph  with ${n\choose 2}-k$ edges remaining virtual.       
The graph is  grown by
selecting an element   uniformly at random   from the pool      of    $N-n+ {n\choose 2}-k$ virtual vertices and edges and changing  status of the sampled element from virtual  to occupied.
Eventually the complete graph on $N$ vertices is built.

We shall focus on the edge-counting  in discrete time, while
identifying  the whole hours $1,\ldots,N$ 
%temporal variable 
with the number of occupied vertices.
The vertices  will be  labelled $1,\ldots,N$  as they get sampled. 
The edges will be directed as $i\leftarrow j$ for $1\leq i<j\leq N$.
We define $G_n$  to be   the  graph on $[n]:=\{1,\dots,n\}$ emerging as the $n$th vertex becomes occupied, when also  $n-1$ virtual edges $i\leftarrow n$ are added to the pool. Thus
at time $n$ the  pool is comprised of some number of unoccupied  edges between vertices in $G_n$, and of $N-n$ virtual vertices.
The events of occupying the edges are viewed as intermediate steps that occur between the whole hours when
 the clock rings and vertices get occupied.
For the main part of the paper we shall consider the process up to the time  $N$, when the last vertex gets occupied and the pool is left with some number of virtual edges.

For the late stages of the process, when a few virtual vertices remain, it is more insightful to switch to a time reversal and vertex-counting.

The analytic tractability of the model relies upon exchangeability intrinsic to the sampling process.
Indeed,
whenever the pool at some stage has $r$ elements, their all $r!$ relative orders in which the elements become occupied are equally likely.
  Let $I_n(i, j)$ denote the indicator of the event that
$G_n$ 
has (occupied) edge $i\leftarrow j$, and let 
\begin{equation}\label{ind-rep}
X_n=\sum_{\{(i,j):~1\leq i<  j<n\}} I_n(i,j), ~~~n=1,\dots,N
\end{equation}
be the number of edges in $G_n$.
We have, by the exchangeability,
\begin{equation}\label{edgeP}
{\mathbb P}[I_n(i, j)=1]=\frac{n-j}{N-j+1},~~~1\leq i<j\leq n.
\end{equation}
The bias in edge probabilities is seen from this formula.
For instance, at times $n$ close to $N/2$ the edge $1\leftarrow 2$ is occupied with probability about $1/2$, while the edge $1\leftarrow n$ is occupied with probability $0$.

The expected value of $X_n$ follows straight from  (\ref{edgeP}):
\begin{eqnarray}\nonumber
\mu_n:={\mathbb E}[X_n]&=&\sum_{j=1}^{n-1} \frac{(j-1)(n-j)}{N-j+1}\\
\label{meanX}
&=& {n\choose 2} +(n-1)(N-n)-N(N-n+1)\sum_{j=1}^{n-1}\frac{1}{N-j+1}\,.
\end{eqnarray}

In what follows we shall work out  formulas and  asymptotics for the variance, the latter being remarkably explicit.
Our central results will capture  more delicate features of the edge-counting process $(X_n)$, for $N$ large.
These include Poisson approximations for  times   close to the start and the end of the process, a fluid limit 
and Gaussian fluctuations for the bulk of the process suitably scaled.
The following diagram illustrates the time scales in the succession of   ${N\choose 2}+N$ occupancy events
that occur as the trivial graph grows to the complete graph on $N$ vertices,
with $\circ$ marking the vertex occupancy and $\bullet$ the edge occupancy events
$$\underbrace{\circ \cdots\circ\bullet  \circ\cdots\circ}_{N^{1/3}-{\rm scale,\, \bullet\, Poisson}}
\underbrace{\bullet 
\bullet\bullet\bullet\cdots\bullet\bullet\bullet\bullet\circ   \bullet\bullet\cdots\bullet\bullet\circ  \bullet   \bullet\bullet  \cdots\bullet\bullet}_{N-{\rm scale,\,\bullet\, Gaussian}}   
\underbrace{
  \circ\bullet  \cdots\bullet \circ   \bullet  \cdots\bullet }_{N\log N-{\rm scale,\,\circ\, Poisson}}$$
In the first two regions the scaling refers to the whole hours marked $\circ$. In the third region 
$\circ$ are rare, occuring at intervals 
 of the order of $N\log N$.

{\it Notation:}  we shall denote $\stackrel{d}{\to}$ the convergence in distribution, $\stackrel{{\mathbb P}}{\to}$ the convergence in probability,
and $\Rightarrow$ the weak convergence in the Skorokhod space of cadlag processes \cite{Billingsley}.

\section{A dynamic urn model}\label{DUM}

The edge-counting can be thought of 
 in  terms of  the following urn process. 
The urn contains white and black balls. A ball is sampled uniformly at random and its colour noticed.
If it was a black ball, the ball is removed without replacement. If it was the $n$th sampled white ball, it is removed and $n-1$ new black balls are replaced to the urn. Assuming that initially the urn had  $N$ white balls and zero black,
 $X_n$ can be identified with the number of black balls removed from the urn before  $n$  white balls are sampled.

It is clear from this representation that $(X_n)$ is a Markov chain.
Let
$$\Delta X_n:=X_n-X_{n-1}$$
denote the increment,
where $X_0=0$. 
Given $X_n=k$,  the increment $\Delta X_{n+1}$ is the number of  black balls  sampled without replacement before a white ball gets drawn
 from an urn containing ${n\choose 2}-k$ black and $N-n$ white balls.
This number has a negative hypergeometric distribution
\begin{equation}\label{TP}
{\mathbb P}[\Delta X_{n+1}= \ell\,|\,X_n=k]=\frac{{M-\ell-1\choose K-1 }}{{M\choose K }},~~~~0\leq \ell\leq M-K,
\end{equation}
with moments
\begin{eqnarray}\label{mean-comp}
{\mathbb E}[\Delta X_{n+1}\,|\,X_n=k]&=&\frac{M-K}{K+1},\\
\label{V}
{\rm Var \,}[\Delta X_{n+1}\,|\,X_n=k]&=&\frac{(M+1)(M-K)K}       {(K+1)^2(K+2)},
\end{eqnarray}
where for shorthand
$$
 M=N-n+{n\choose 2} -k,~~~K=N-n.
$$

We note  a particularly simple formula for the second difference
\begin{equation}\label{delta-di}
\me[\Delta_{n+1}-\Delta_n]=\mu_{n+1}-2\mu_n+\mu_{n-1}=\frac{n-1}{N-n+1}.
\end{equation}
This can be neatly explained in terms of  a decomposition of the edge-counting process.
For $m<n$, the pool at time $n$ contains some of the {\it aged} virtual edges missing in the graph $G_m$ and some of the {\it recent}  virtual edges added to the pool between times $m+1$ and $n$.
The difference $X_n-X_m$ decomposes accordingly as
\begin{equation}\label{OR}
X_n-X_m=A_{m,n}+R_{m,n}
\end{equation}
where $A_{m,n}$ is the number of aged edges   occupied by time $n$.
The increments 
$$A_{m,m+1}, A_{m,m+2}-A_{m,m+1},\ldots, A_{m,n}-A_{m,n-1}$$
are exchangeable, with  equal mean values.
In particular, for $m=n-1$, the mean contribution to $\Delta_{n}$ and $\Delta_{n+1}$ of the aged  edges is the same, so the increase in (\ref{delta-di}) occurs due to some of $n-1$
 edges added to the pool at time $n$ and occupied by time $n+1$.

\section{Two martingales}

Let $({\cal F}_n^X)$ be the natural filtration associated with $(X_n)$.  By the Markov property,  (\ref{mean-comp}) is equivalent to
\begin{equation}\label{mean-comp-1}
{\mathbb E}\left[\Delta X_{n+1}|{\cal F}_n^X\right]=\frac{{n\choose 2}-X_n}{N-n+1}.
\end{equation}
The predictable process
\begin{equation}\label{compensator}
C_n:=\sum_{j=1}^{n-1} \frac{{j\choose 2}-X_j}{N-j+1}
\end{equation}
is the  compensator involved in the Doob-Meyer decomposition of $(X_n)$, i.e. the sequence $(X_n-C_n)$ is a martingale.

Taking expectation  in  (\ref{mean-comp-1}) yields a recursion for the mean 
\begin{equation}\label{rec-mu}
\mu_{n+1}-\mu_n=\frac{{n\choose 2}-\mu_n}{N-n+1},
\end{equation}
which could have been concluded directly from (\ref{meanX}).
Introduce the centered process 
$$X_n^\circ:=X_n-\mu_n.$$
Subtracting (\ref{rec-mu})  from (\ref{mean-comp-1})  we arrive at
\begin{equation}\label{sec-mart}
{\mathbb E}\left[\frac{X_{n+1}^\circ}{N-n}{\bigg|}{\cal F}_n^X\right]=\frac{ X_{n}^\circ}{N-n+1},
\end{equation}
which says that the sequence $(N-n+1)^{-1}X_n^\circ,~n=1,\ldots,N,$ is also a martingale.

\section{Asymptotics of the moments}

To approximate the mean  (\ref{meanX}) for  large $N$ consider  the function
$$
\varphi(t):=2\int_0^t \frac{s(t-s)}{1-s}\,{\rm d}s=  2(1-t)L+2t-t^2,
$$
where and henceforth
$$L:=\log(1-t).$$
A counterpart of (\ref{rec-mu}) is the differential equation
$$\varphi'(t)=\frac{t^2-\varphi(t)}{1-t},~~~~~\varphi(0)=0.$$

From the second expression in (\ref{meanX}) one readily derives the estimate
\begin{equation}\label{mean-a}
\max_{1\leq n\leq N} \left|\frac{2}{N^2}\,\,\mu_n-\varphi\left( \frac{n}{N}\right)\right| = O\left( \frac{1}{N}\log\left(\frac{N}{N-n+1} \right)\right).
\end{equation}
whichis uniform in $n$. 
This entails the convergence rate $O(N^{-1})$  for $n/N$ bounded away from $1$, that is 
in the range $n/N<1-\delta$ for any given $\delta\in (0,1)$.

We turn next to the variance
$$\sigma_n^2:={\rm Var}[X_n].$$
To derive a recursion, condition on ${\cal F}_n^X$ to decompose $\sigma_{n+1}^2$ as
\begin{eqnarray*}
\sigma_{n+1}^2&=& {\mathbb E} \{{\rm Var}[X_{n+1}|{\cal F}_n^X]\}+{\rm Var}\{{\mathbb E}[X_{n+1}|{\cal F}_n^X]\}\\
&=&
{\mathbb E} \{{\rm Var}[\Delta X_{n+1}|{\cal F}_n^X]\}+{\rm Var}\left[\frac{N-n}{N-n+1}X_{n}\right],
\end{eqnarray*}
where the second expression follows from (\ref{sec-mart}). Now, from (\ref{V})

\begin{eqnarray*}
{\mathbb E} \{{\rm Var}[\Delta X_{n+1}|{\cal F}_n^X]\}=\\
 \frac{N-n}{(N- n+1)^2(N-n+2)}&&\!\!\!\!\!\!\left[ \sigma_n^2 +\left({n\choose 2}-\mu_n\right)\left(N-n+1+{n\choose 2}-\mu_n \right)\right],
\end{eqnarray*}
and the second term is 
$$
{\rm Var}\left[\frac{N-n}{N-n+1}X_{n}\right]=\left(\frac{N-n}{N-n+1}\right)^2  \sigma_n^2.$$
Adding up and appealing to (\ref{rec-mu}) with a bit of algebra we arrive at the recursion 
\begin{equation}\label{rec-sigma}
\sigma_{n+1}^2-\sigma_n^2=\frac{(N-n)(\mu_{n+1}-\mu_n)(\mu_{n+1}-\mu_n+1)-2\sigma_n^2}{N-n+2}.
\end{equation}
Solving (\ref{rec-sigma}) gives
\begin{equation}\label{sigma-sum}
\sigma_n^2 = (N-n+1)(N-n+2) \sum_{j=1}^{n-1} \frac{(\mu_{j+1}-\mu_j)(\mu_{j+1}-\mu_j+1) }{(N-j+1)(N-j+2)}.
\end{equation}

From (\ref{meanX}) the difference of the mean values is
\begin{equation}\label{diff-means}
\mu_{n+1}-\mu_n=\frac{-(n-1)(N-n)}{N-n+1}  +N\sum_{i=1}^{n-1} \frac{1}{N-i+1}.
\end{equation}
Approximating the sum  in (\ref{diff-means}) by an integral gives an upper bound
$$\mu_{n+1}-\mu_n<N \log\left(\frac{N}{N-n+1} \right),$$
which substituted in   (\ref{sigma-sum}) readily gives $\sigma_n^2=O(N^3)$ uniformly in $n$.  The latter does not capture the decay of the variance for $n/N$ close to $1$.
With some  more work for the upper half of the range we obtained the estimate
\begin{equation}\label{var-est}
\sigma_n^2 < 5(N-n+2)N^2 \log^2 \left(\frac{N}{N-n+1}\right)~~{\rm for~}~n>N/2,
\end{equation}
which shows the right order of the variance (though the constant $5$ is not optimal).

A continuous counterpart of (\ref{rec-sigma}) is the differential equation
$$
\psi'(t)=-\frac{2\,\psi(t)}{1-t}+\frac{1}{4}(\varphi'(t))^2,~~~\psi(0)=0,
$$
whose solution is 
\begin{equation}\label{psi}
\psi(t)=
(1-t)\{(2-t)L^2+2(3-t)L+t(6-t)\}.
\end{equation}
This approximates the variance in the sense that
$$
\max_{1\leq n\leq N} |N^{-3}    \sigma_n^2-\psi(n/N)|\to  0,{\rm ~~~as~}N\to\infty,
$$
where the approximation error is $O(N^{-1})$ for $n/N<1-\delta$. 
We have $\psi(1)=0$ which neatly agrees with the $n=N$ instance of (\ref{var-est}).

We note   that the analogue of (\ref{sigma-sum}) is the integral representation
\begin{equation}\label{psi-int}
\psi(t)=(1-t)^2\int_0^t \left(\frac{s+\log(1-s)}{1-s} \right)^2{\rm d}s.
\end{equation}

\begin{figure}

\includegraphics[scale=0.64]{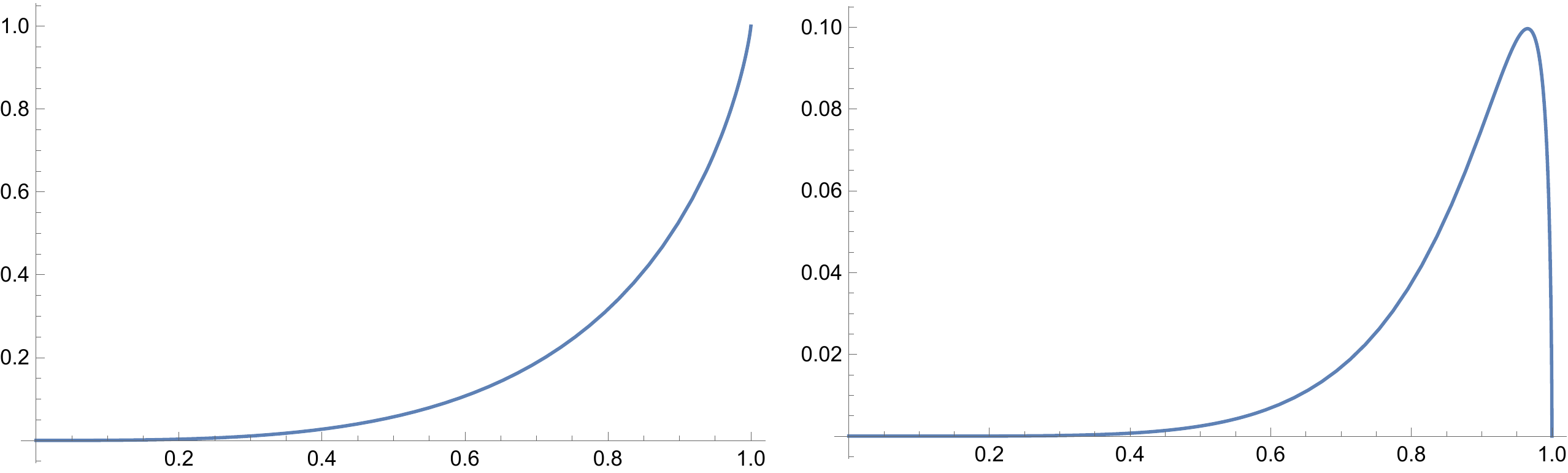}

\caption{$\varphi$ and $\psi$. Observe different scales!}
\end{figure}

It is instructive to arrive at the asymptotic variance (\ref{psi}) in a more direct way,  in order to emphasise the  dependence among the edges in $G_n$.
 Consider $a\leftarrow i$ and $b\leftarrow j$ with $i<j<n$.  
The edge  $a\leftarrow i$ is not occupied at time $n$ with probability $(N-n+1)/(N-i+1)$.
Given this event, the probability that the edge $b\leftarrow j$ is not occupied at time $n$ equals  $(N-n+2)/(N-j+2)$,
which is lower than the unconditional probability.
With this in mind, the second moment of the number of virtual edges ${n\choose 2} -X_n$ 
is given by
\begin{eqnarray}
\nonumber
{\mathbb E}
\left[{n\choose 2} -X_n\right]^2&=&{n\choose 2} -\mu_n\\
\nonumber
&+&(N-n+1)(N-n+2)    \sum_{i=1}^n \frac{(i-1)(i-2)}{(N-i+1)(N-i+2)} \\
\label{secmom}
&+&2(N-n+1)(N-n+2) 
\sum_{1\leq i<j\leq n}  \frac{(i-1)(j-1)}{(N-i+1)(N-j+2)}
\end{eqnarray}
and subtracting the squared mean
\begin{eqnarray*}\left({\mathbb E}\left[{n\choose 2} -X_n \right]  \right)^2&=& \\
(N-n+1)^2\sum_{i=1}^n\left(\frac{i-1}{N-i+1} \right)^2   &+&2 (N-n+1)^2 \sum_{1\leq i<j\leq n}  \frac{(i-1)(j-1)}{(N-i+1)(N-j+1)}  
\end{eqnarray*}
the variance 
$$\sigma_n^2={\rm Var}\left[{n\choose 2} -X_n\right]$$
is given by
\begin{eqnarray*}\sigma_n^2=
{n\choose 2}-\mu_n-\sum_{i=1}^n (i-1)\left(\frac{N-n+1}{N-i+1} \right)^2+(N-n+1)\sum_{i=1}^n \frac{(i-1)(i-2)(n-i)}{(N-i+2)(N-i+1)^2}
 \\
2 (N-n+1) \sum_{1\leq i<j\leq n}  \frac{(i-1)(j-1)(n-j)}{(N-i+1)(N-j+1)(N-j+2)}  
\end{eqnarray*}
For $N\to\infty$ the last term with a bivariate sum is of the dominating order $O(N^3)$. Letting $N\to\infty$, the scaled term becomes
\begin{eqnarray}\nonumber
\psi(t)=2(1-t)\iint\limits_{0<x<y<t}\frac{xy(t-y)}{(1-x)(1-y)^2}{\rm d}x{\rm d}y,
\end{eqnarray}
which is yet another integral representation for $\psi$.
We see that the dominating contribution to the variability of $X_n$ is due to the edges with different origins.

On summary, we obtained the asymptotics
\begin{equation}\label{asymp-mom}
\mu_n\sim \frac{N^2}{2} \varphi(n/N),~~~\sigma_n^2\sim N^3\psi(n/N),
\end{equation}
valid uniformly for $n/N<1-\delta$.

\section{Poisson approximation for the early stages} \label{Early}

For small $n$ the graph $G_n$ is likely to contain no edges.
From the first expression in (\ref{meanX}) it is seen that
 \begin{equation}\label{Omean}
\mu_n=O\left(\frac{n^3}{N} \right)
\end{equation}
uniformly in $n\leq N$, and that the order estimate is sharp.  Borrowing the terminology from \cite{RG}, we can say  that $N^{1/3}$ is a threshold function  for the existence of edges in $G_n$.
To obtain a more precise result, 
let the random variable
$$
\xi_N :=\min\{n: X_n>0\}
$$
be the first time when an edge becomes occupied. That is to say, the graph $G_n$ for $n<\xi_N$  has no edges, while $G_{\xi_N}$ has exactly one.

\begin{theorem}\label{P1} As $N\to\infty$ 
$$N^{-1/3}\,\xi_N\stackrel{d}{\to} \xi,$$
where ${\mathbb P}[\xi>x]=e^{-x^3/6},~~~x>0$. Thus $\xi^3$ has exponential distribution with mean $6$.
The convergence  holds with all moments, that is for $\alpha\geq 0$
$$
{\mathbb E}[N^{-\alpha/3}\xi_N^\alpha]\to {\mathbb E}[\xi^\alpha]=6^{\alpha/3}\,\Gamma\left(1+{\alpha}/{3}\right).
$$
\end{theorem}

\proof
Using (\ref{TP})
$${\mathbb P}[\xi_N>n]=\prod_{j=2}^{n-1}  \frac{N-j}{{j\choose 2}+N-j}=\exp\left\{    -\sum_{j=2}^{n-1} \log\left(1+\frac{{j\choose 2}}{N-j} \right)\right\} $$

%For $N$ large,  expanding the logarithm we have uniformly in $j<n<N^{2/5}$
%$$\log\left(1+\frac{{j\choose 2}}{N-j} \right)=\frac{{j\choose 2}}{N-j}+O\left(\frac{j^4}{N^2} \right).$$
Splitting
$$\frac{1}{N-j}=\frac{1}{N}+\frac{j}{N(N-j)}$$
we have in the range $j<n<N^{2/5}$
$$\frac{{j\choose 2}}{N-j}=\frac{{j\choose 2}}{N}+O\left(\frac{j^3}{N^2} \right), $$
and so
$$\log\left(1+\frac{{j\choose 2}}{N-j} \right)=\frac{{j\choose 2}}{N}+O\left(\frac{j^4}{N^2} \right), $$
whence summing up 
$$\sum_{j=2}^{n-1} \log\left(1+\frac{{j\choose 2}}{N-j} \right)=\frac{n^3}{6N}+O\left(\frac{n^2}{N} \right) + O\left(\frac{n^5}{N^2} \right) $$
uniformly in 
 $n<N^{2/5}$. 
Letting $n\sim x N^{1/3}$ this converges to $x^3/6$, thus
$${\mathbb P}[N^{-1/3}\xi_N>x]\to e^{-x^3/6}$$
and the convergence in distribution follows.

To prove  convergence of the moments we need to check the uniform integrability condition 
\begin{equation}\label{UI}
\limsup\limits_{N\to\infty} {\mathbb E}[N^{-\alpha/3}\xi_N^\alpha ~1(N^{-1/3} \xi_N >A)]\to 0 ~~~{\rm as~}A\to\infty.
\end{equation}
The  expectation involved is
\begin{eqnarray}\label{Tr}
N^{-\alpha/3}\sum_{AN^{1/3}<n\leq N} n^\alpha ~{\mathbb P}[\xi_N=n],
\end{eqnarray}
where 
$${\mathbb P}[\xi_N=n+1]=        \frac{{n\choose 2}}{{n\choose 2}+N-n} ~  \prod_{j=2}^{n-1}  \frac{N-j}{{j\choose 2}+N-j}$$

For $n>N^{1/2}$ we estimate
\begin{eqnarray*}
{\mathbb P}[\xi_N=n+1]<{\mathbb P}[\xi_N>n]<{\mathbb P}[\xi_N>N^{1/2}]=\\
\exp\left\{ -\sum_{1\leq j<  N^{1/2}}   \log\left(1+\frac{{j\choose 2}}{N-j} \right)\right\}<e^{-c N^{1/2}}
\end{eqnarray*}
for some $c>0$. Splitting the sum (\ref{Tr}) at $N^{1/2}$ for the upper part we have

$$\sum_{N^{1/2}\leq n\leq N}   
 n^\alpha ~{\mathbb P}[\xi_N=n]
< \sum_{ N^{1/2}\leq n\leq N}  n^\alpha e^{-c n^{1/2}}<N^{\alpha+1}e^{-c N^{1/2}}\to 0.$$
For the lower part, using the bounds
$$ \frac{{j\choose 2}}{N-j}<2,~~~~~\log(1+x)>\frac{x}{2}~~{\rm for~~}0<x<2,$$
we estimate for large enough $N$
\begin{eqnarray*}
\sum_{AN^{1/3}< n< N^{1/2}}   
 n^\alpha ~{\mathbb P}[\xi_N=n]&<&\sum_{AN^{1/3}< n< N^{1/2}} \frac{n^{\alpha+2}}{2N}\exp\left\{- \sum_{j=1}^{n-1}\log\left(     1+
\frac{{j\choose 2}}{N-j}
\right)\right\}<\\
\sum_{AN^{1/3}< n< N^{1/2}} \frac{n^{\alpha+2}}{2N}
\exp\left\{- \sum_{j=1}^{n-1}   
\frac{{j\choose 2}}{2N}\right\}&=&
\sum_{AN^{1/3}< n< N^{1/2}} \frac{n^{\alpha+2}}{2N}
\exp\left\{-    
\frac{n(n-1)(n-2)}{6N}\right\}<\\
\int\limits_{AN^{1/3}}^\infty   \frac{x^{\alpha+2}}{N} ~ \exp\left(-\frac {x^3}{6N}\right)\,{\rm d}x&=&N^{\alpha/3}\int\limits_{A^3}^\infty   \frac{y^{\alpha/3}}{3} ~ \exp\left(-\frac {y}{6}\right)\,{\rm d}y.
\end{eqnarray*}
Multiplying by $N^{-\alpha/3}$ and sending $A\to\infty$,  condition (\ref{UI}) follows.
\qed

More informative is the following Poisson approximation on the $N^{1/3}$ scale, which actually implies the convergence in distribution in Theorem \ref{P1}.

\begin{theorem}\label{P2} As $N\to\infty$,
the process $(X_{\lfloor t N^{1/3}\rfloor},~t\geq 0)$ converges weakly in the Skorohod space $D[0,\infty)$ to a nonhomogeneous Poisson process with rate function $\lambda(t)=t^2/2, ~t>0$.
\end{theorem}
\proof  We shall verify two conditions of Theorem 5 from \cite{GS}, Ch. 5, Section 3.

Fix $T>0$. For each $i$, the probability that by time $TN^{1/3}$ there is an occupied edge $\cdot\leftarrow i$ does not exceed 
$T^2 N^{-1/3}$,
and the probability of at  least two occupied  edges (possibly with same origins) is estimated as the square of this quantity. Summing over all pairs $\cdot\leftarrow i, \cdot\leftarrow j$ in the bounds $$sN^{1/3}<i,j<tN^{1/3}$$
 we obtain
$${\mathbb P}[X_{\lfloor t N^{1/3}\rfloor}-X_{\lfloor s N^{1/3}\rfloor}\geq 2]=O((t-s)^2), ~~~0<s<t<T.$$
That is to say, asymptotically the time-scaled process has no jumps bigger than $1$. This is the first condition of the cited theorem.

The compensator of $(X_n)$  is given by (\ref{compensator}).
For $n\leq tN^{1/3}$ and $N$ large we have the expansion
$$C_n=\sum_{m=1}^{n-1} \frac{{m\choose 2}}{N}-\sum_{m=1}^{n-1} \frac{X_m}{N-m+1}+O\left(N^{-1/3}\right),$$
where the second sum converges in probability to $0$, since
from (\ref{Omean})
$${\mathbb E}\left[\sum_{1\leq m<tN^{1/3}} \frac{X_m}{N-m+1}\right]< \sum_{1\leq m<tN^{1/3}} \frac{c \,m^3}{N(N-m+1)}=O(N^{-2/3}).$$
Hence as $N\to\infty$ for each $t$
$$C_{\lfloor tN^{1/3}\rfloor}\stackrel{{\mathbb P}}{\to}\frac{t^3}{6}.$$
Thus the time-scaled compensator converges to a deterministic function. This is the second condition  of  the cited theorem. 
%By monotonicity of the compensator it follows from  this implies that for every $T>0$
%$$\sup_{t\in[0,T]}\left| C_{tN^{1/3}}-\frac{t^3}{6}    \right| \stackrel{{\mathbb P}}{\to} 0.$$
\qed

\section{The fluid limit}

Fluid limits are functional laws of large numbers which assert 
 convergence of a process to a deterministic curve \cite{Darling}.
Comparison with the mean readily implies a fluid limit for the scaled edge count. 
\begin{theorem} As $N\to\infty$,
\begin{equation}\label{fl-lim}
\sup_{t\in[0,1]}\left| \frac{2}{N^2}\,X_{\lfloor tN\rfloor}-\varphi(t)\right|  \stackrel{{\mathbb P}}{\to} 0.
\end{equation}
\end{theorem}

\proof From the asymptotics (\ref{asymp-mom})
$${\rm Var}\left[\frac{2}{N^2}\,X_{\lfloor tN\rfloor}\right]\to 0,$$
hence  application  of Chebyshev's inequality yields
$$\frac{2}{N^2}\,X_{\lfloor tN\rfloor}\stackrel{{\mathbb P}}{\to } \varphi(t),$$
for each fixed $t\in (0,1)$. But this implies the functional convergence
 in the uniform metric, because
$X_{n}$ is nondecreasing in $n$ and the function $ \varphi(t)$  is continuous (and monotone). 
The underlying fact  is that the pointwise convergence of nondecreasing functions on $[0,1]$ to a continuous function
ensures the uniform convergence.
\qed

\section{Gaussian fluctuations}

The typical magnitude of the difference in (\ref{fl-lim}) is of the order of $N^{-1/2}$. To capture such fluctuations we need to scale with factor $N^{1/2}$,
with an eye on 
the Central Limit Theorem
$$\frac{X_{\lfloor tN\rfloor}-\frac{N^2}{2}\varphi(t)}{\sqrt{\psi(t) N^3}}\stackrel{d}{\to} {\cal N}(0,1),$$
for each fixed $t\in(0,1)$. 
This indeed holds true, but is not straightforward since the indicators in the representation (\ref{ind-rep}) are not independent.

A wider view is to  consider the whole process
\begin{eqnarray}\nonumber
Y_N(t):&=&\frac{\sqrt{N}}{2}\left(\frac{2}{N^2}X_{\lfloor tN\rfloor}-\varphi(t)\right)\\
&=& \label{YN}
\frac{X_{\lfloor tN\rfloor}-\frac{N^2}{2}\varphi(t)}{N^{3/2}},~~~t\in[0,1]
\end{eqnarray}
as a random element of the Skorohod space $D[0,1]$ of piecewise continuous functions.
We aim to show weak  convergence to a Gaussian diffusion  defined by the stochastic integral
$$Y(t)=(1-t)\int_0^t \frac{s+\log(1-s)}{1-s}{\rm\, d}B(s),~~~t\in[0,1],$$
where $B(\cdot)$ is the Brownian motion. The process $Y(\cdot)$ can be characterised as the unique (strong) solution to the
 stochastic differential equation
\begin{equation}\label{SDE}
{\rm d}Y(t)= -\frac{Y(t)}{1-t}{\rm d}t+ (L+t){\rm d}B(t),~~~Y(0)=0.
\end{equation}
At the right  endpoint the value is nonrandom, $Y(1)=0$.
This process  is a sort of Gaussian bridge derived from suitably time-changed Brownian motion, as in  \cite{Gasbarra} (see second construction in Example 1).

Using It{\^o} isometry the covariance function is computed as
\begin{eqnarray}\nonumber
{\rm Cov}(Y(s), Y(t))&=& \frac{1-t}{1-s}  \int_0^s\left( \frac{u+\log(1-u)}{1-u}\right)^2\,{\rm d}u\\
\label{cov}
&=& \frac{1-t}{1-s}\,\psi(s),~~~0\leq s\leq t\leq 1,
\end{eqnarray}
which agrees and generalises the asymptotic variance formula (\ref{psi-int}) obtained differently.
Thus we may also  characterise $Y(\cdot)$ as a continuous Gaussian process on $[0,1]$ with zero mean and the covariance function (\ref{cov}).

\begin{theorem} As $N\to\infty$ 
$$(Y_N(t),~t\in[0,1])\Rightarrow (Y(t),~t\in[0,1]),$$ 
where the weak convergence of processes holds in the Skorohod space $D[0,1]$. 
\end{theorem}

\proof  Replacing centering in (\ref{YN}) by the mean $\mu_n$ makes no difference asymptotically, by the virtue of (\ref{mean-a}).
%$$\frac{2}{N^2}\mu_{tN}-\varphi\left(t\right)=O(N^{-1}).$$
Thus consider the naturally centered process
$$Y^\circ_N(t):=\frac{X_{\lfloor tN\rfloor }-\mu_{\lfloor tN\rfloor }}{N^{3/2}}=    \frac{X_{\lfloor tN\rfloor }^\circ}{N^{3/2}}.$$
The generator ${\cal L}_N$ of ${Y}_N^\circ (\cdot)$ acts on smooth functions   as
\begin{equation}\label{generator}
{\cal L}_N f(t,x)=N\,{\mathbb E}[f(t+N^{-1}, Y^\circ_N(t+N^{-1}))-f(t,Y^\circ_N(t))\,|\,Y^\circ_N(t)=x].
\end{equation}

Setting $t=n/N$, $\Delta X_{n}^\circ=X_{n}^\circ-X_{n-1}^\circ$ we have
\begin{equation}\label{incr}
Y^\circ_N\left(t+N^{-1}\right)-Y^\circ_N(t)= \frac{\Delta X_{n+1}^\circ }{N^{3/2}}.
\end{equation}
Given 
$$X_n=\mu_n+ x N^{3/2} ~~~~(\Longleftrightarrow Y^\circ_N(t)=x)$$
for some admissible $x$, the distribution of increment $\Delta X_{n+1}$ is negative hypergeometric  with moments (\ref{mean-comp}), (\ref{V}).
Recalling (\ref{rec-mu}),
%$$\frac{{n\choose 2}-\mu_n}{N-n+1}-(\mu_{n+1}-\mu_n)=0,$$
the conditional expectation  of the increment (\ref{incr}) expands as
\begin{equation}\label{first}
{\mathbb E}\left[ \frac{\Delta X_{n+1}^\circ}{N^{3/2}}\,{\bigg |}\,X_n=\mu_n+x N^{3/2} \right]=\frac{-x}{N-n+1}=\frac{-x}{1-t}N^{-1}+ o\left(N^{-1}\right).
\end{equation}
This quantity squared is of the  order $O(N^{-2})$, hence the second conditional moment of the increment is approximable by the  conditional variance as

$${\mathbb E}\left[ \left(\frac{\Delta X_{n+1}^\circ }{N^{3/2}}\right)^2     \,{\bigg |}\,X_n=\mu_n+x N^{3/2} \right                          ]=
{\rm Var}\left[ \frac{\Delta X_{n+1}}{N^{3/2}}\,{\bigg |}\,X_n=\mu_n+x N^{3/2} \right]+o\left(N^{-1}\right).
$$
Manipulating (\ref{V}), the asymptotics of  the conditional variance becomes 
\begin{equation}\label{second}
\frac{\left({n\choose 2} -\mu_n\right)^2}{(N-n)^2 N^3}+o(N^{-1})=\frac{(t^2-\varphi(t))^2}{4(1-t)^2N}+o(N^{-1})=\frac{(t+L)^2}{N}+o(N^{-1}).
\end{equation}
Expanding the function in (\ref{generator}) in Taylor series,
and applying (\ref{first}) and (\ref{second}),  we may write  the generator (\ref{generator}) in the form
$${\cal L}_N f(t,x)= {\cal L} f(t,x)+o(1),   ~~N\to\infty,$$ 
where ${\cal L}$ is the generator of the diffusion  (\ref{SDE}) acting as
$${\cal L} f(t,x)=
f_t(t,x)+\frac{-x}{1-t}  f_x(t,x)+  \frac{1}{2}(t+L)^2 f_{xx}(t,x).$$

By some general theory (\cite{Kallenberg}, Proposition 17.9 and Theorem 17.25), convergence of generators on a rich enough class of test  functions  implies the functional convergence of processes.
This principle allows us to justify the weak convergence $Y_N^\circ\Rightarrow Y$ 
for the processes restricted to  $[0,\,1-\delta], 0<\delta<1$.
Indeed, within the reduced time range
 the drift and diffusion coefficients in (\ref{SDE}) are Lipschitz-continuous, hence  by \cite{Kolokoltsov} (Theorem 3.17)
 the space of  $C^3$-smooth functions  vanishing with all derivatives at infinity 
 is an invariant  core for the the semigroup associated with (\ref{SDE}).

At the right endpoint $t=1$ the drift and  diffusion coefficients both  explode. To circumvent this nuisance 
 we shall use the fact that  the limit process has the fixed right endpoint $Y(1)=0$, hence 
by the continuity   
for all positive $x,\varepsilon$ there exists $\delta>0$ such that 
\begin{equation}\label{little}
\prob \left[\sup_{1-\delta\leq t\leq 1} |Y(t)|>x\right]<\varepsilon
\end{equation}
(as is concluded from the tightness criterion in \cite{Billingsley}, Theorem 8.2).
Obviously, the inequality 
 $\sup_{1-\delta\leq t\leq 1} |Y(t)|<x$ implies that the modulus of continuity of $Y$ on $[1-\delta,1]$ does not exceed $x$. 
Therefore, with the weak convergence in $D[0,\,1-\delta]$ in hand,
tightness of the sequence $Y_N^\circ$ in $D[0,1]$ will be shown  if we could verify  a counterpart of (\ref{little}) for $Y_N^\circ$ in order to bound the modulus of continuity,
in line with condition (15.9) from \cite{Billingsley}.

For the latter, we shall adopt the Chow  maximal inequality \cite{Gut} (Theorem 9.3):
for positive submartingale $(Z_n)$ and nonincreasing sequence $(c_n)$ it holds that for $\lambda> 0$
$$
\prob\left[ \max_{1\leq k\leq n} c_k Z_k>\lambda \right]\leq \frac{1}{\lambda}\left(\sum_{k=1}^{n-1} (c_k-c_{k+1}) \me[ Z_k] +c_n\me [Z_n]\right).
$$
Applying the inequality to the squared martingale (\ref{sec-mart}),
$$Z_n=\left( \frac{X_n^\circ }{N-n+1}\right)^2,$$
with constants $c_n=(N-n+1)^2$ and $\lambda= x^2 N^3$ we obtain 

\begin{eqnarray*}
\prob\left( \max_{N(1-\delta)\leq n\leq N}    |X_n^\circ|>x N^{3/2} \right) \leq \frac{1}{x^2 N^3} \sum_{n=\lfloor N(1-\delta)\rfloor}^N  \frac{2N-2n+1}{(N-n+1)^2}\,\sigma_n^2   <\\
\frac{4}{x^2 N} \sum_{n=\lfloor N(1-\delta)\rfloor}^N   \log^2\left(\frac{N}{N-n+1}\right)\to  \frac{4}{x^2 }\, \delta \log^2\delta,
\end{eqnarray*}
as $N\to\infty$.
Choosing $\delta$ small enough, this implies the analogue of (\ref{little}) for $Y_N^\circ$, hence the tightness in $D[0,1]$.
Since  weak convergence in each $D[0,\,1-\delta]$ entails the convergence of finite-dimensional distributions, 
the assertion now follows by application of Prohorov's theorem \cite{Billingsley}.\qed

A similar 
 technical difficulty caused by discontinuity of the drift coefficient at $t=1$  is present also by a martingale approach to approximating the empirical  distributions by the Brownian bridge.
In \cite{Jacod} (p. 561) this is handled using the time reversibility, which is not available in our situation.

\section{Construction by insertion}

\paragraph{The P{\'o}lya urn.}
For $M$ balls,  of which $K$ are white and $M-K$ are black, arrangement in a sequence defines a weak composition of  integer $M-K$ in $K+1$ parts. 
The white balls appear in the role of delimeters between the parts.
For instance
\begin{equation}\nonumber
\bullet\bullet\bullet  \stackrel{1}{\circ} \bullet\bullet\bullet\,\bullet \stackrel{2} {\circ} \,  \stackrel{3} {\circ} \bullet\bullet\bullet\bullet\bullet   \stackrel{4}{\circ}\bullet 
\end{equation}
represents the weak composition $(3,4,0,5,1)$ of integer $13=3+4+0+5+1$.
We speak of `weak' composition to mean that zero parts are allowed.

Suppose $K$ white balls 
%labelled $1,\ldots, K$ 
are arranged in sequence and a black ball is inserted uniformly at random in any of $K+1$ possible positions, that is to the left, to the right or between two balls of the existing configuration. Then another black ball is inserted uniformly at random
in any of $K+2$ possible positions.  Continuing so forth yields a sequence of random weak compositions, each  with 
uniform distribution over ${M\choose K}$ weak compositions in $K$ parts for $M=K, K+1,\ldots$.
Each part of the composition has a negative binomial distribution as in (\ref{TP}), and as $M$ grows the whole composition behaves like a $(K+1)$-component P{\'o}lya urn process
(see \cite{Borovkov} for a recent study).
In particular, as $M\to\infty$, the vector of parts scaled
by $M-K$ converges (almost surely) to a vector $(\xi_1,\ldots,\xi_{K+1})$ with uniform distribution over the $K$-dimensional simplex.

In application to our random graph, the P{\'o}lya process is directly  relevant to the analysis  of the aged component $A_{m,n}$ in (\ref{OR}), and this connection will be employed in the next Section.
To represent the full edge counts  $(X_n)$ a more sophisticated insertion procedure is needed.

\paragraph{The insertion algorithm}\label{Insert}
Start with $N$ white balls arranged in a row and labelled $1,\ldots,N$. 
In hours $1$ and $2$ do nothing.
In hour $3$ insert a black ball
uniformly at random in any position to the right of  the $2$nd white ball. 
By induction, within hour $n=3,\ldots,N$ do  one-by-one $n-2$ such random insertions of  black balls in the configuration to the right of $(n-1)$st white ball.
For instance, if $N=6$ the configuration 
\begin{equation}\nonumber \label{balls}
 \stackrel{1}{\circ} \,\stackrel{2} {\circ} \,  \stackrel{3} {\circ} \bullet\,\bullet   \stackrel{4}{\circ}\bullet\,\bullet  \stackrel{5}{\circ}\bullet \stackrel{6}{\circ}
\end{equation}
may occur in hour $5$, when a black balls remains to be inserted to the right of the $4$th while ball to complete the step, and there are $6$ equiprobable ways to do this.

\vskip0.2cm

In this picture, the composition emerging by the end of hour $N$ has the same distribution as 
$(\Delta X_1,\ldots,\Delta X_N)$. 
In terms of the graph process, the difference with the urn model in Section \ref{DUM} is that the virtual elements become scheduled 
for future occupancy as soon as they are added to the pool. 
Implicitly, the insertion algorithm has been used in our computation of the moments of $X_n$.

\section{Poisson approximation for the last stages}

Finally, we seek for approximation of the edge-counting process when a few  vertices remain unoccupied.
Complementing our notation, let
$$\Delta X_{N+1}:= {N\choose 2}-X_N$$
be the number of edges missing in $G_N$. Interpreting $\Delta X_{N+1}$ as the number of edges  occupied after the last vertex becomes occupied,
this quantity is analogous to the number $\xi_N$
%=\min\{n: X_n>0\}$ 
of vertices occupied before the first edge becomes occupied.

To find the right scaling we shall inspect the moments of $\Delta X_{N+1}$.
For the mean we have from (\ref{meanX}) a simple formula

$$  {\mathbb E}[\Delta X_{N+1}]={N\choose 2}-\mu_N=N\sum_{j=1}^{N-1}\frac{1}{N-j+1}=N(h_N-1),$$
where $h_N$ is the harmonic number. Thus as $N\to\infty$
\begin{equation}\label{last-me}
  {\mathbb E}[\Delta X_{N+1}]\sim N\log N.
\end{equation}
For the second noncentral moment (\ref{secmom}) gives
\begin{eqnarray*} {\mathbb E}[\Delta X_{N+1}]^2= {\mathbb E}[\Delta X_{N+1}]+ 2\sum_{i=1}^N \frac{(i-1)(i-2)}{(N-i+1)(N-i+2)}+\\
2\sum_{1\leq i<j\leq N} \frac{2(i-1)(j-1)}{(N-i+1)(N-j+2)}.
\end{eqnarray*}
Approximating the first sum by an integral gives
$$2\sum_{i=1}^N \frac{(i-1)(i-2)}{(N-i+1)(N-i+2)}=O(N^2).$$
We are left with assessing
$$2\sum_{1\leq i<j\leq N} \frac{2(i-1)(j-1)}{(N-i+1)(N-j+2)}.$$
Replacing in the denominator $1$ by $2$ will give a lower bound, and $2$ by $1$ an upper bound; so let us stick with the denominator $(N-i+1)(N-j+1)$. 
Then adding the same sum over $i>j$ and the diagonal part
$$
2\sum_{i=1}^N \frac{2(i-1)^2}{(N-i+1)^2}=O(N^2),
$$
the lower bound becomes
$$\sum_{i=1}^N\sum_{j=1}^N  \frac{2(i-1)(j-1)}{(N-i+1)(N-j+1)}\sim  2\left(N\sum_{j=1}^{N-1}\frac{1}{N-j+1}\right)^2\sim 2 (N\log N)^2.$$
The upper bound has the same asymptotics, whence
$$  {\mathbb E}[\Delta X_{N+1}]^2\sim 2(N\log N)^2.$$
By the same token, for the third moment the leading term will be
$${\mathbb E}[\Delta X_{N+1}]^3\sim 3! \sum_{1\leq i<j<k\leq N} \frac{3!(i-1)(j-1)(k-1)}{(N-i+1)(N-j+2)(N-k+3)}\sim 3!(N\log N)^3,$$
and proceeding so forth we obtain
$${\mathbb E}[\Delta X_{N+1}]^r\sim r!(N\log N)^r.$$
Since $r!$ is the $r$th moment of the exponential distribution, we conclude that
$$\frac{\Delta X_{N+1}}{N\log N}\stackrel{d}{\to}{\rm Exponential}(1).$$

Comparing the mean (\ref{last-me}) with the number  $N-1$ of edges added to the pool at time $N$, we see that the latter has no impact on the limit law of $\Delta X_{N+1}$.
Similarly, for each fixed $m$,  in the decomposition (\ref{OR}) of $X_{N}-X_{N-m}$ the term $R_{N-m,N}$ counting the recent edges is asymptotically negligible. This idea 
combined with exchangeability 
leads to the following result.

%Note that  ${N\choose 2}-X_N$ is the number of edges missing in $G_N$. 
%Looking at the terminal configuration of white and black balls for $N$ large,  black balls occur at the left edge  by a nonhomogeneous Poisson process, with scaling the position number by factor $N^{-1/3}$. 
%We show next that on the right edge scaling by $(N\log N)^{-1}$ is required, leading to occurrence of white balls by a homogeneous Poisson process. 

\begin{theorem}\label{conv-last} For each $m\geq 0$ the random vector
\begin{eqnarray*}
\label{prop-vec}
\left(\frac{\Delta X_{N-m+1}}{N\log N},\, \frac{\Delta X_{N-m}} {N\log N},\, \ldots, \,\frac{\Delta X_{N+1}} {N\log N}\right)
\end{eqnarray*}
converges in distribution, as $N\to\infty$, to a vector of $m+1$ i.i.d.  exponential random variables with mean one.
\end{theorem}
%The argument has two ingredients.
%Firstly, using the method of moments we show that, with this scaling, the limit law of ${N\choose 2}-X_{N-m}$ is Gamma$(m+1,1)$. Secondly a splitting result will be used, 
%very much in spirit of the P{\'o}lya urn limit at start of the Section.
\proof
For $m>0$ we consider  
$${N\choose 2} - X_{N-m}=  \Delta X_{N-m+1}+\cdots+\Delta X_{N+1},$$
the number of edges in the complete graph minus the number of egdes in $G_{N-m}$. This includes
${N-m\choose 2} -X_{N-m}$ aged edges  
and 
$${N\choose 2}-{N-m+1\choose 2} = (m+1)N-{m+2\choose 2}=O(N)$$
recent edges, all relative to time $N-m$.
From (\ref{meanX})
$${\mathbb E} \left[{N-m\choose 2} -X_{N-m}\right]\sim (m+1)N\log N,$$
hence with scaling by $N\log N$ we may ignore the recent ages.

From (\ref{secmom}) for the second moment we get
\begin{eqnarray*}
{\mathbb E} \left[{N-m\choose 2} -X_{N-m}\right]^2\sim
2 (m+1)(m+2) 
\sum_{1\leq i<j\leq N-m}  \frac{(i-1)(j-1)}{(N-i+1)(N-j+2)}\sim\\
 (m+1)(m+2) (N\log N)^2.
\end{eqnarray*}
It is now easy to recognise and prove the asymptotics of higher moments:
$${\mathbb E} \left[{N-m\choose 2} -X_{N-m}\right]^r\sim  (m+1)\cdots (m+r) (N\log N)^r.$$
Recalling that $(m+1)_r=(m+1)\cdots(m+r)$ are moments of the gamma  distribution, we  
conclude that
%$$\frac{{N-m\choose 2} - X_{N-m}}{N\log N}\stackrel{d}{\to} {\rm Gamma}(m+1),$$
%and so
$$\frac{{N\choose 2} - X_{N-m}}{N\log N}\stackrel{d}{\to} {\rm Gamma}(m+1,1).$$
This special case is also known as the Erlang distribution.

Appealing to the insertion construction in Section \ref{Insert}, given that
${N\choose 2} - X_{N-m}=M-m$ virtual edges are present at time $N-m$, the weak composition 
corresponding to their occupancy schedule
can be represented by the
$(m+1)$-component P{\'o}lya urn with $M-m$ black balls. Since for $N\to \infty$ also ${N\choose 2} - X_{N-m}\stackrel{\prob}{\to} \infty$,
the proportions converge in distribution to random vector $(\xi_1,\ldots,\xi_{m+1})$ uniformly distributed over the $m$-dimensional simplex.
It follows that
  the vector (\ref{prop-vec}) converges in distribution
%$$\frac{\Delta X_{N-m+1}}{N\log N}, \frac{\Delta X_{N-m}} {N\log N}, \ldots, \frac{\Delta X_{N+1}} {N\log N}$$
to
$(\gamma\xi_1,\ldots,\gamma\xi_{m+1})$,
where $\gamma\stackrel{d}{=}{\rm Gamma}(m+1,1)$
is independent of $(\xi_1,\ldots,\xi_{m+1})$.
By Lucacs  characterisation of the gamma distribution, such $(\gamma\xi_1,\ldots,\gamma\xi_{m+1})$ has i.i.d. Exponential$(1)$ components.
\qed

As an easy corollary of this result, for the time reversal of $(X_n)$ we obtain a counterpart of Theorem \ref{P2}:

\begin{theorem} As $N\to\infty$,
the point process with atoms
$$\frac{{N\choose 2}-X_N}{N\log N}, \,\frac{{N\choose 2}-X_{N-1}}{N\log N},\,\cdots$$
converges weakly  to a homogeneous Poisson process with unit rate.
\end{theorem}

In particular, for $x>0$, at the stage when the complete graph on $N$ points falls short of  $\lfloor x N\log N\rfloor$ elements (vertices and edges),
the distribution of the number of virtual vertices is close to Poisson$(x)$.

%https://doi.org/10.1007/978-3-540-70847-6_15

\end{document}